\documentclass[graybox]{svmult}

\usepackage{type1cm}        
\usepackage{makeidx}         
\usepackage{graphicx}        
\usepackage{multicol}        
\usepackage[bottom]{footmisc}

\usepackage{newtxtext}       %
\usepackage[varvw]{newtxmath}       

\makeindex             

\begin{document}

\title*{Hierarchical Coarse Basis by Randomised SVD: the Helmholtz Problem}
\author{Martin J. Gander\orcidID{0000-0001-8450-9223} and\\ Yao-Lin
  Jiang\orcidID{0000-0002-5541-1136} and\\ Hui Zhang\orcidID{0000-0001-7245-0674}}
\institute{Martin J. Gander \at Department of Mathematics, University of Geneva, Rue du
  Conseil-Général 9, CP 64, 1211 Genève 4, Switzerland, \email{martin.gander@unige.ch} \and Yao-Lin
  Jiang \at School of Mathematics and Statistics, Xi'an Jiaotong University, 710049 Xi'an, China,
  \email{yljiang@mail.xjtu.edu.cn} \and Hui Zhang (corresponding author) \at School of Mathematics
  \& Physics, Xi'an Jiaotong-Liverpool University, Ren'ai Road 111, 215123 Suzhou, China,
  \email{hui.zhang@xjtlu.edu.cn}}
%
%
\maketitle

\abstract*{The oscillatory waves require sufficient degrees of freedom to resolve. That restriction
  usually applies also to coarse problems for Schwarz methods. The resulting coarse problem is then
  too large. To address the issue, a new form of Schwarz methods with coarse correction is proposed
  for the Helmholtz problem. There are two components in the proposed form: randomised SVD of
  interface iteration, and hierarchical domain decomposition. The resulting coarse problem is
  hierarchical and can be solved in a decoupled way at each level of hierarchy.}

\section{Introduction}
\label{sec:1}

The Helmholtz problem is hard to solve by established iterative methods for positive definite
problems. Special techniques have been developed to make domain decomposition methods (DDM) work for
the Helmholtz problem. However, the efficiency of DDMs is still restricted by the oscillatory nature
of propagating waves. One big concern is the large size of coarse problems for coupling many small
subdomains. This motivated us to rethink the traditional framework of coarse correction which aims
at small subdomains and treats large coarse problems recursively by DDM. An alternative way is to
start with big subdomains and thus a small coarse problem, then treat big subdomains recursively by
DDM. This leads to a hierarchy of subdomains and coarse problems.

\begin{figure}
  \centering
  \includegraphics[scale=.43,trim=0 0 0 0,clip]{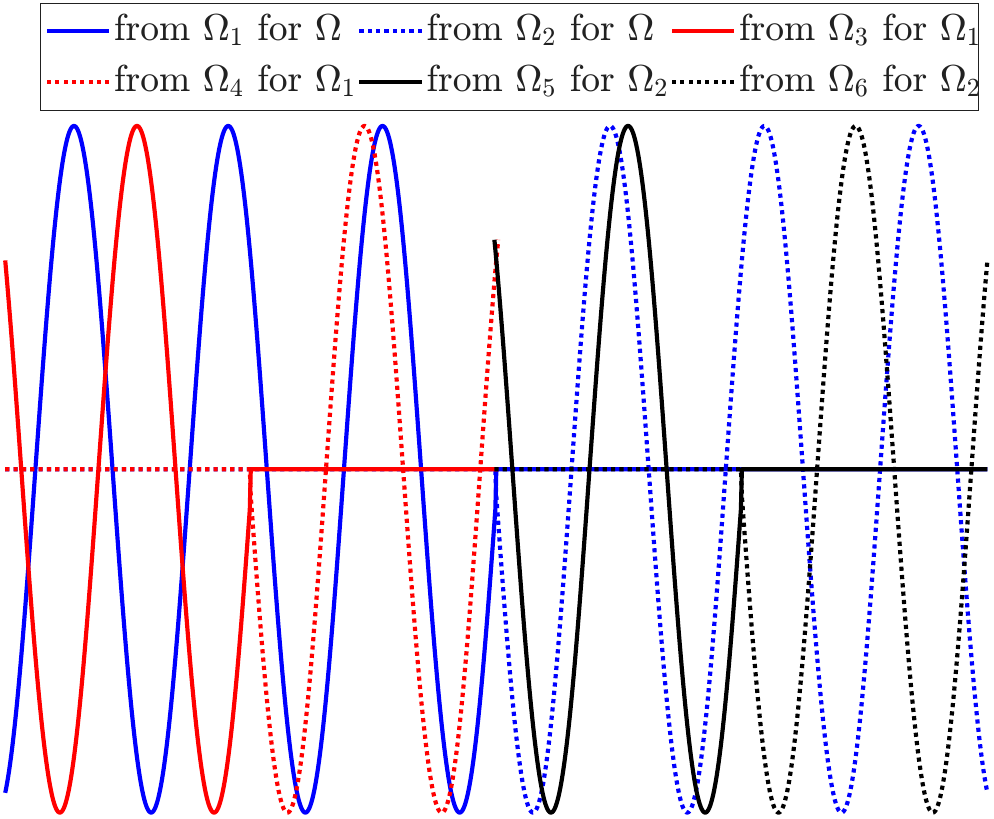}
  \caption{Hierarchical coarse basis $\{C,C_1,C_2\}$ (real parts) from 2-level bisections \&
    restricted $P_i$.}
  \label{fig1}
\end{figure}

The hierarchical coarse basis is generated very naturally by this approach. To see this, we consider
the Helmholtz problem in 1D with impedance boundary conditions:
\begin{equation}\label{lap1}
  -u''-k^2u=f\text{ in }\Omega=(0,1),\quad -(u'+\mathrm{i}ku)(0)=g_0,\;(u'-\mathrm{i}ku)(1)=g_1.
\end{equation}
We first decompose $\Omega$ into $\Omega_{1}=(0,b_{1})$ and $\Omega_{2}=(a_{2},1)$ with
$a_{2}<b_{1}$. Let $a_{1}=0$ and $b_2=1$. Suppose (\ref{lap1}) after discretisation becomes
$A\mathbf{u}=\mathbf{f}$. Then we can solve $A\mathbf{u}=\mathbf{f}$ by the optimised Schwarz
iteration with coarse correction
\begin{eqnarray}\label{ras1}
  \tilde{\mathbf{u}}^{(n)}&=&(P_1\tilde{A}_1^{-1}R_1+P_2\tilde{A}_2^{-1}R_2)(\mathbf{f}-A\mathbf{u}^{(n-1)}),\label{sub1}\\
  \mathbf{u}^{(n)}&=&\tilde{\mathbf{u}}^{(n)}+CA_c^{-1}C^H(\mathbf{f}-A\tilde{\mathbf{u}}^{(n)}),\label{coarse1}
\end{eqnarray}
where $R_i$ is the identity restriction from $\overline{\Omega}$ to $\overline{\Omega_i}$, $P_i$ is
a prolongation from $\overline{\Omega_i}$ to $\overline{\Omega}$ such that $P_1R_1+P_2R_2$ is the
identity, $\tilde{A}_i$ is from discretisation of the Helmholtz equation in $\Omega_i$ with
impedance boundary conditions $-(u'+\mathrm{i}ku)(a_i)=0$, $(u'-\mathrm{i}ku)(b_i)=0$, $A_c=C^HAC$
is the Galerkin approximation of $A$ in the range of $C$ (the column space of $C$), and
$C=[{\color{red}C_1^0,C_2^0}]$ with ${\color{red}C_i^0}=P_iU_i$ with {\color{red}$U_i\approx u_i$}
obtained from solving the discretised version of
{\color{red}\begin{eqnarray}
  -u_1''-k^2u_1&=&0\text{ in }\Omega_1,\quad -(u_1'+\mathrm{i}ku_1)(0)=0,\;(u_1'-\mathrm{i}ku_1)(b_1)=1,\label{U1}\\
  -u_2''-k^2u_2&=&0\text{ in }\Omega_2,\quad -(u_2'+\mathrm{i}ku_2)(a_2)=1,\;(u_2'-\mathrm{i}ku_2)(1)=0,\label{U2}
\end{eqnarray}}
which is well-known as the $a$-harmonic\footnote{Here, $a$-harmonic means the weak form zero source
  problem with the bilinear form $a(u,v)$.} extension of the interface data. The motivation of using
$a$-harmonic extensions is that the error of the optimised general additive Schwarz iteration
(\ref{sub1}) for $n\ge 1$ is essentially coming from incorrect interface values at $b_1$ and $a_2$;
see \cite{oras}. In other words, the error space is exactly the $a$-harmonic extension of the
interface space. In particular for (\ref{sub1}), the error space is of dimension two, and the coarse
correction (\ref{coarse1}) would result in the exact solution of (\ref{lap1}). To get the
hierarchical domain decomposition method, we now perturb (\ref{sub1}) and (\ref{U1})-(\ref{U2}) by
solving the problems on $\Omega_1$ and $\Omega_2$, again with the method
(\ref{sub1})-(\ref{coarse1}). For example, we decompose $\Omega_1$ into $\Omega_3=(0,b_3)$ and
$\Omega_4=(a_4,b_1)$, and construct $C_1=[{\color{red}C_3^1,C_4^1}]$. If the problems on $\Omega_3$,
$\Omega_4$ are solved exactly, then the method of (\ref{sub1})-(\ref{coarse1}) applied to the
$\Omega_1$ problem would converge in one step. Or, we can further decompose $\Omega_3$, $\Omega_4$,
and so on. Note that the set of coarse basis $\{C,C_1,C_2\}$ or $\{C,C_1,C_2,C_3,C_4,C_5,C_6\}$
(with $\Omega_i$ decomposed into $\Omega_{2i+1}$ and $\Omega_{2i+2}$, and $C_i$ the basis on
$\Omega_i$) after prolongation forms a hierarchical coarse basis on $\Omega$. An example of the
hierarchical Helmholtz coarse basis is shown in Fig.~\ref{fig1}, where we used the restricted
prolongation $P_i$ (see \cite{ras, oras}) from optimised restriced additive Schwarz (ORAS).

\section{Hierarchical Coarse Basis in 2D}
\label{sec:2}
We now extend the idea introduced in Sec.~\ref{sec:1} to 2D. Since there are many dofs on an
interface of a 2D domain, it is no longer efficient to lift the whole interface space to the coarse
space. The key question is what kind of interface modes to pick. We follow the general idea of using
the error modes that are difficult to contract by the Schwarz iteration; see \cite{eigenSchwarz,
  eigenRobinSchwarz}.

\begin{figure}[t]
  \includegraphics[scale=.32,trim=0 0 0 0,clip]{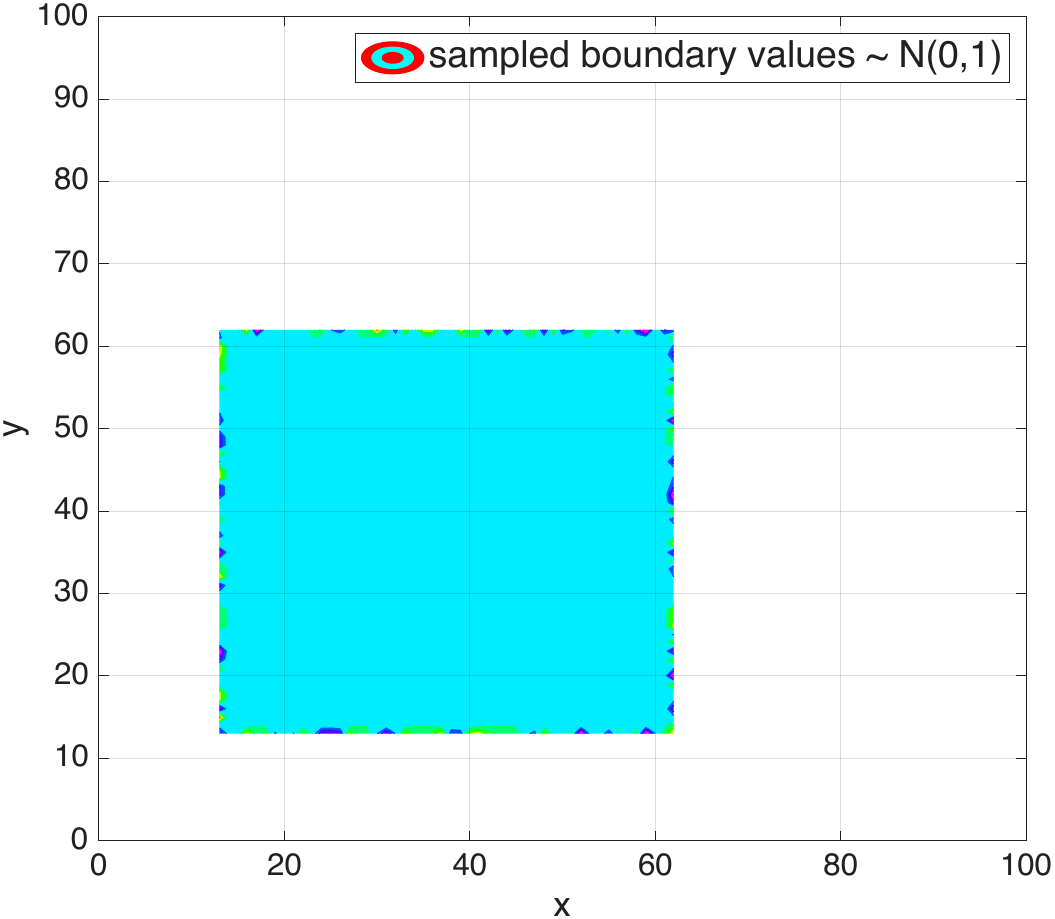}
  \includegraphics[scale=.34,trim=32 8.51 12 0,clip]{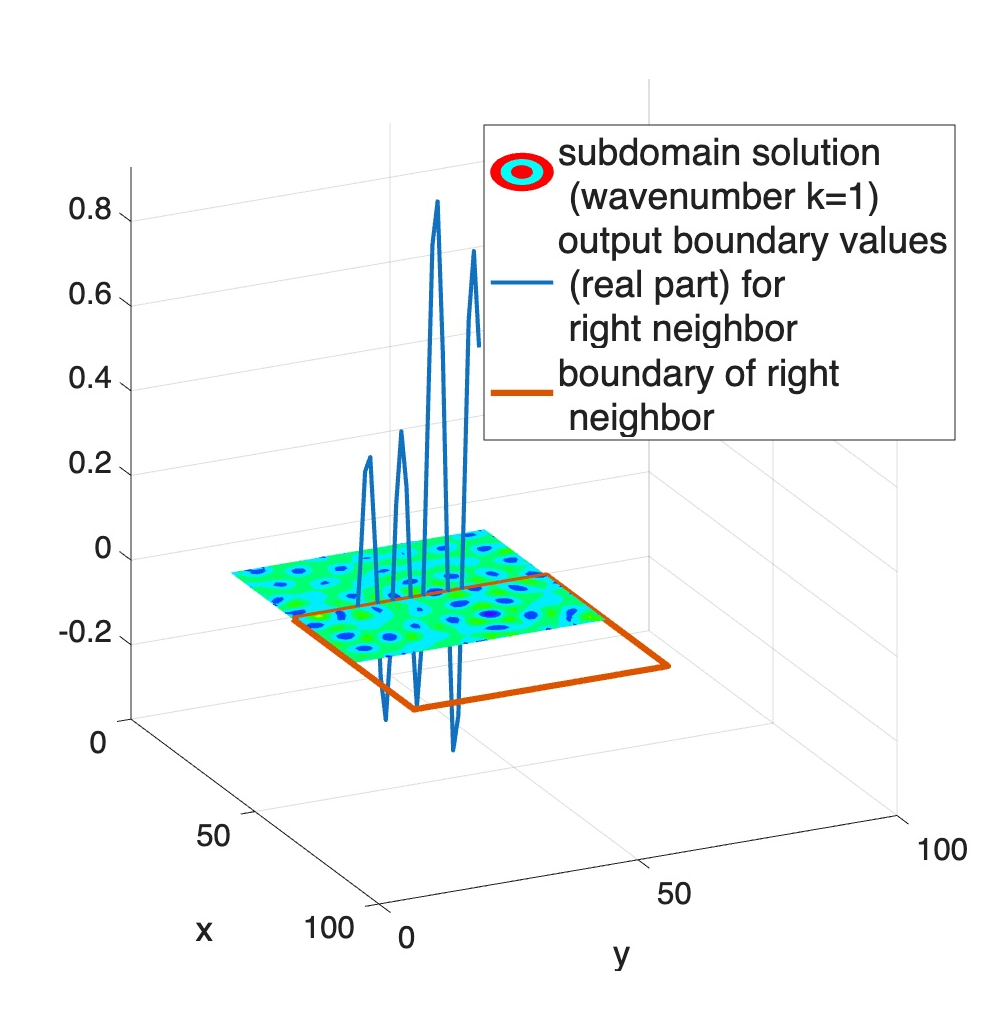}
  \caption{Interface map from one subdomain to its right neighbor with boundary operator
    $\partial_{\mathbf{n}}-\mathrm{i}k$.}
  \label{fig2}
\end{figure}

To this end, for each subdomain $\Omega_i$ (or each interface
$\overline{\partial\Omega_i\cap\Omega_j}$), we consider the interface map $T_i$ (or $T_{ij}$) as
composition of: 1) the $a$-harmonic extension of the interface data $g_i$ on
$\overline{\partial\Omega_i\cap\Omega}$ (or $g_{ij}$ on $\overline{\partial\Omega_i\cap\Omega_j}$)
to the subdomain solution $u_{i}$, followed by 2) the evaluation of interface data $\{g_{ji}\}$ on
$\{\overline{\partial\Omega_j\cap\Omega_i}\}$ for the neighboring subdomains $\{\Omega_j\}$ using
$u_i$. The interface map $T_i$ (or $T_{ij}$) is a ``rectangular'' operator in the sense that its
domain and codomain are different spaces. An example of the interface map from one subdomain to one
of its neighbors is shown in Fig.~\ref{fig2}.

By \cite{despres}, $\sum_{j\ne i}\|g_{ji}\|_{L^2}^2$ should be a fraction of $\|g_i\|_{L^2}^2$ (or
$\|g_{ij}\|_{L^2}^2$), which gives contraction. To find out the slowly convergent error modes, we
use the Singular Value Decomposition\footnote{Ma-Alber-Scheichl \cite{ma2023} used the SVD of a
  restriction operator to build a coarse basis. {\color{red}Similar works are \cite{hu2024novel,
      nataf2024coarse, parolin2025coarse}.}} e.g. $T_iV=U\Sigma$. Then the columns of $V$
corresponding to large singular values (diagonal entries of $\Sigma$) are the coarse modes on the
interface of $\Omega_i$. Since the SVD is very expensive and we need only a few dominant singular
vectors, we shall use the randomised SVD (rsvd) instead; see \cite{rsvd10, rand2011}, which needs
only to apply $T_i$ on a few random samples (a bit more than the required number of modes),
orthogonalise, and apply the adjoint $T_i^*$ before the SVD.

{\color{red}In this work, we adopt $T_i$ in a purly algebraic form without considering function spaces, and
simply use the Euclidean vector inner product. The Q1 finite element method is used for
discretisation. Let $B_{i}$ be the boolean matrix that restricts the nodal values $\mathbf{u}_i$ on
$\overline{\Omega}_i$ to the boundary $\partial\Omega_i$. Then the discretised Robin boundary value
problem with zero source on $\Omega_i$ is in the form $\tilde{A}_i\mathbf{u}_i=B_i^T\mathbf{g}_i$
where $\mathbf{g}_i$ comes from the Robin boundary value. Denote the solution by
$\mathbf{u}_i^{*}$. Let $B_{ij}$ be the matrix that maps $\mathbf{u}_i^{*}$ to the Robin boundary
data $\mathbf{g}_j$ for $\Omega_j$. A convenient definition is $B_{ij}=B_j\tilde{A}_jR_jP_i$ with
sufficient overlap and $P_i$ so that the Robin trace (residual) on
$\partial\Omega_j\cap\overline{\Omega}_{i}$ could be evaluated from the side of
$\overline{\Omega}_j\cap\overline{\Omega}_i$ rather than $\overline{\Omega}_i-{\Omega}_j$. We note
that $B_{ij}$ has many zero rows corresponding to the boundary values on
$\partial\Omega_j-(\partial{\Omega}_j\cap\overline{\Omega}_i)$, and the definition of $B_{ij}$ could
be refined to not take the zero rows. Let $B_{i}^o=[B_{ij}]$ be the vertical concatenation for all
the neighbors of $\Omega_i$. The matrix $T_i$ is defined as $T_i=B_i^o\tilde{A}_i^{-1}B_i^T$. After
finding the dominant right singular vectors $V_i$ of $T_i$, we extend them to the coarse vectors
$C_i=P_i\tilde{A}_i^{-1}B_i^TV_i$. The coarse basis matrix is defined as the horizontal
concatenation of $C_i$'s for all the subdomains, namely $C=(C_i)$.
}
\begin{remark}
  It is well known that the propagating modes need to travel through the outgoing boundaries of the
  domain to disappear. A consequence is that if there are $m$ subdomains in one direction across the
  domain to the outgoing boundaries, then it takes $m$ parallel subdomain iterations to see a
  significant decrease of the error. So the interface map defined through the input and output of
  \emph{multiple} Schwarz iterations could be useful for extracting coarse modes, which is still to
  be investigated.
\end{remark}

{\color{red}\begin{remark}
  In the literature \cite{ma2023,hu2024novel,nataf2024coarse,parolin2025coarse} for the Helmholtz
  problem, some weighted $H^1$-like norms are used for the $a$-harmonic extension of the boundary
  values. That could be useful in our setting too. We observed for the Laplace problem that some
  singular corner modes can ``pollute'' the desired smoothing modes (e.g. the piecewise constant
  modes), and that can be removed by using $H^{1/2}$-norm instead of the $L^2$-norm.
\end{remark}
}
It is worthwhile to mention that in a hierarchical decomposition, e.g. dividing each of the
$2$-by-$2$ subdomains into $2$-by-$2$ smaller subdomains, the coarse modes are built with a few
steps of the sublevel iteration like (\ref{sub1})-(\ref{coarse1}) for the $a$-harmonic extension,
and thus contain iterative errors.



\begin{figure}[t]
  \sidecaption[t]
  \includegraphics[scale=.35,trim=0 0 0 0,clip]{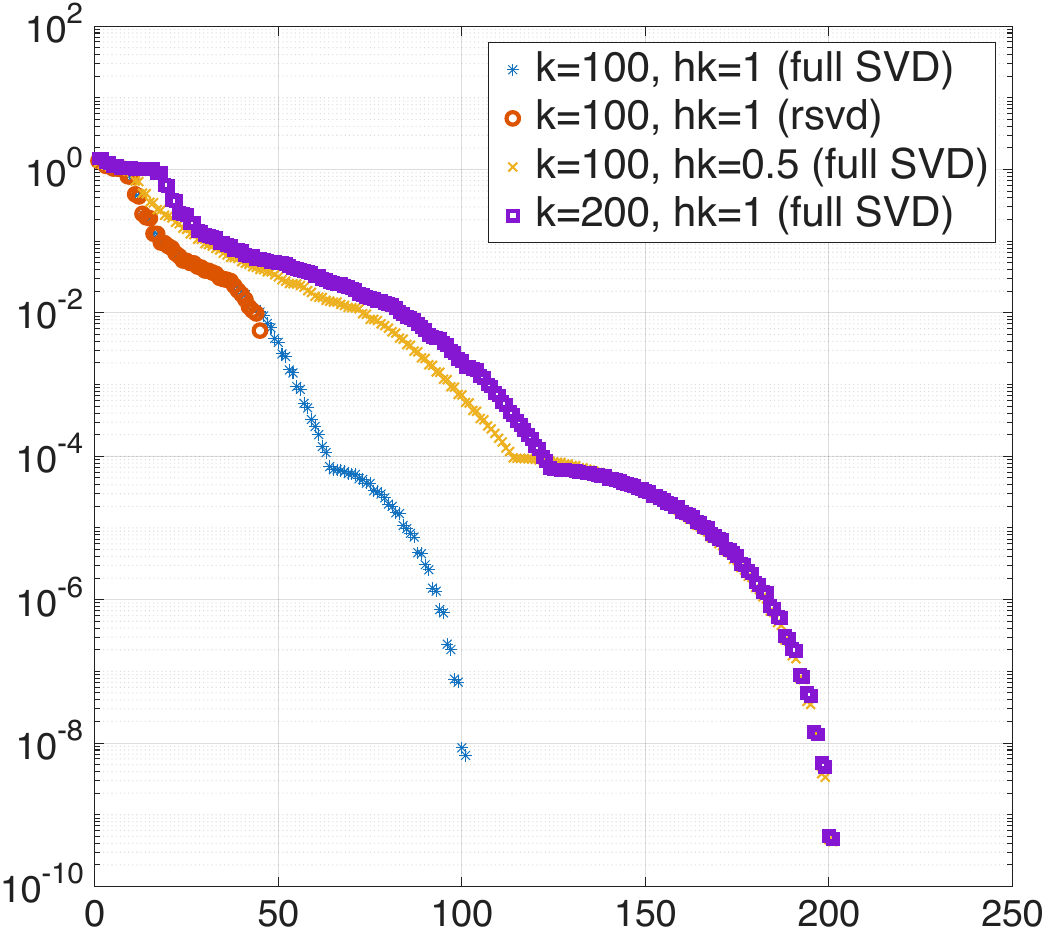}
  \caption{Singular values of the interface map $T_i$ on the meshes of element size $h$ from one of
    the $2\times 2$ subdomains to its neighbors with $4h$ overlap, for $\Delta+k^2$ in
    $\Omega=(0,1)^2$ with impedance conditions. The horizontal axis is the index of the descending
    singular values. The actual values of the singular values are along the vertical axis.}
  \label{fig4}
\end{figure}

The singular values of the interface map from the subdomain to its neighbors are plotted in
Fig.~\ref{fig4}. We see the initial singular values decay very slowly, which can be attributed to
propagating modes. Their number grows linearly with the wavenumber $k$. This is a big difference
from the diffusion problem. There are also two regimes the decaying is slowing down in the
middle. Otherwise, the singular values decay exponentially, just as for the diffusion problem, which
can be attributed to the evanescent modes. The number of singular values above a fixed value grows
about linearly with wavenumber $k$ as $hk=1$, and hence linearly with subdomain size $H$ by a
scaling argument. This number is also inversely proportional to $h$ (essentially the overlap tending
to zero). In this sense, it is useful to keep the overlap large and independent of $h$. However, the
price is solving larger subdomain problems.

\section{Numerical Experiments}
\label{sec:3}
We test the hierarchical restricted additive Schwarz method with the coarse basis generated by the
randomised SVD, as described in Sec.~\ref{sec:2}.  We denote the hierarchical {\em decomposition} by
$[m_{x1}\times m_{y1}, m_{x2}\times m_{y2}, \ldots,m_{x\ell}\times m_{y\ell}]$ meaning the domain is
first partitioned into $m_{x1}\times m_{y1}$ subdomains, then each subdomain is partitioned into
$m_{x2}\times m_{y2}$ subdomains, and so on. For convenience, we use ${m}$ without subscript when
$m_{x1}=m_{y1}=\ldots$. The {\em number of coarse modes} from a subdomain at each level of partition
is collected as $[n_{c1},n_{c2},\ldots,n_{c\ell}]$, which means e.g.  ${n_{c1}}$ coarse modes are
generated for each subdomain of the partition $m_{x1}\times m_{y1}$. We take $n_{ci}+5$ samples in
rsvd. The {\em number of Schwarz iterations} like (\ref{sub1})-(\ref{coarse1}) used at level $l$ for
both ``the coarse basis at the upper level'' and ``the preconditioner for the original problem'' is
denoted by ${n_{il}}$. For convenience, we use $n_{c}$ and $n_{i}$ if the quantities are level
independent. A bilinear finite element method on rectangular meshes is used for the
discretisation. The partition is element based and first without {overlap}, then extended outward by
{two-element} width from the non-overlapping interfaces, i.e. {\em overlap} ${\delta=4h}$ with $h$
the element size. Each non-overlapping subdomain at the last level $\ell$ is discretised with
$n\times n$ elements, i.e. ${H/h=n}$ with $H$ the subdomain size. We use GMRES \cite{gmres} to solve
the non-symmetric indefinite Schwarz preconditioned system. {\color{red}In all the experiments,
  $k_{\max}h=1$ is kept, where $k_{\max}$ is the maximum wavenumber in the domain.}

\subsection{Free Space Problem}
\label{sec:3.1}
We solve $-\Delta u-k^2u=f$ in $\Omega=(0,1)^2$, $\partial_{\mathbf{n}}u-\mathrm{i}ku=0$ on
$\partial\Omega$ with a random solution (by using a random $\mathbf{u}$ in the discretised problem
$A\mathbf{u}=\mathbf{f}$ to get $\mathbf{f}$). The initial guess is zero. The stopping criterion is
the relative residual $\|\mathbf{f}-A\mathbf{u}^{(n)}\|_2/\|\mathbf{f}\|_2<10^{-5}$.  The number of
coarse modes per subdomain $n_c$ listed is the least to reach the corresponding iteration number.

We first test with only one level of partition, which amounts to an optimised two-level Schwarz
method with our rsvd coarse basis. The results are shown in Tab.~\ref{tab1}. Note that double $m$ or
$n$ leads to half element size and in turn half overlap and double wavenumber (because $kh$ is
fixed). To keep roughly equal iteration numbers, as we increase the wavenumber $k$ and refine the
mesh such that $kh=1$ on fixed subdomains (double $n$ and fixed $m$), the number of coarse modes per
subdomain $n_c$ needs to grow linearly; but if additionally we increase proportionally also the
number of subdomains $m$ in each direction (double $m$ and fixed $n$), then $n_c$ seems to be stable
for $m\ge 4$. To check finer subdomains for fixed wavenumber and meshes, we look along the skew
diagonals (half $n$ and double $m$), and we find the required $n_c$ roughly proportional to the
subdomain size $H=1/m$ for $m\ge 4$. This is consistent with the observation made from
Fig.~\ref{fig4} in Sec.~\ref{sec:2}.

Some differences from the diffusion problem are: 1) insufficient number of coarse modes for the
Helmholtz problem causes slower convergence than with no coarse modes, and 2) the sufficient number
of coarse modes per subdomain, to be better than without coarse modes, or to keep the same iteration
number, grows linearly with $k$ for fixed $kh$ and subdomains. Actually, to make an improvement over
without coarse modes, $n_c$ needs to exceed the number of propagating modes per subdomain and cover
some evanescent modes. This is much like the rule of thumb for grid based discretisation e.g. ten
points per wavelength, even though the coarse problem here resembles a spectral element method.

\begin{table}[t]
  \caption{GMRES \#iter with $n_{c}$ (below \#iter) coarse modes per subdomain, $n_{i}\!=\!1$
    iteration for the preconditioner, $m^2$ subdomains, $(mn)^2{\color{red}=k^2}$ elements and $4h$
    overlap for the {\bf free space} problem.}
  \label{tab1}
  \begin{tabular}{p{0.4cm}p{2.8cm}p{2.8cm}p{3cm}p{2.3cm}}
    \hline\noalign{\smallskip}
    $n$  & $m=2$ & $m=4$ & $m=8$ & $m=16$ \\    
    \noalign{\smallskip}\svhline\noalign{\smallskip}
    \phantom{1}4 & \phantom{1}6 \quad \phantom{1}5 \quad \phantom{1}4 \quad \phantom{1}3 & 15 \quad \phantom{1}8 \quad \phantom{1}4 \quad \phantom{1}3 & 32 \quad \phantom{1}6 \quad \phantom{1}4 \quad \phantom{1}3 & 61 \quad \phantom{1}8 \quad \phantom{1}4 \quad \phantom{1}3\\
         & \phantom{1}0 \quad \phantom{1}2 \quad \phantom{1}3 \quad \phantom{1}4 & \phantom{1}0 \quad \phantom{1}4 \quad \phantom{1}5 \quad \phantom{1}6 & \phantom{1}0 \quad \phantom{1}6 \quad \phantom{1}7 \quad \phantom{1}8 & \phantom{1}0 \quad \phantom{1}6 \quad \phantom{1}7 \quad \phantom{1}8\\
    \phantom{1}8 & \phantom{1}7 \quad \phantom{1}9 \quad \phantom{1}5 \quad \phantom{1}3 & 17 \quad \phantom{1}8 \quad \phantom{1}5 \quad \phantom{1}3 & 34 \quad \phantom{1}6 \quad \phantom{1}4 \quad \phantom{1}3 & 63 \quad \phantom{1}7 \quad \phantom{1}4 \quad \phantom{1}3\\
         & \phantom{1}0 \quad \phantom{1}1 \quad\phantom{1}5 \quad \phantom{1}7 & \phantom{1}0 \quad 10 \quad 11 \quad 13 & \phantom{1}0 \quad 11 \quad 13 \quad 14 & \phantom{1}0 \quad 11 \quad 13 \quad 14\\
    16  & \phantom{1}9 \quad \phantom{1}7 \quad \phantom{1}6 \quad \phantom{1}3 &  20 \quad \phantom{1}6 \quad \phantom{1}4 \quad \phantom{1}3 & 39 \quad \phantom{1}6 \quad \phantom{1}4 \quad \phantom{1}3 & 78 \quad \phantom{1}6 \quad \phantom{1}4 \quad \phantom{1}3  \\
         & \phantom{1}0 \quad 11 \quad 12 \quad 13 & \phantom{1}0 \quad 19 \quad 23 \quad 25 & \phantom{1}0 \quad 21 \quad 23 \quad 25 & \phantom{1}0 \quad 22 \quad 24 \quad 25\\
    32  &  10 \quad \phantom{1}5 \quad \phantom{1}4 \quad \phantom{1}3 &  22 \quad \phantom{1}6 \quad \phantom{1}4 \quad \phantom{1}3 & 48 \quad \phantom{1}6 \quad \phantom{1}4 \quad \phantom{1}3 & 92 \quad \phantom{1}6 \quad \phantom{1}4 \quad \phantom{1}3 \\
         & \phantom{1}0 \quad 27 \quad 28 \quad 29 & \phantom{1}0 \quad 39 \quad 43 \quad 47 & \phantom{1}0 \quad 42 \quad 45 \quad 48 & \phantom{1}0 \quad 43 \quad 45 \quad 49\\
    \noalign{\smallskip}\hline\noalign{\smallskip}
  \end{tabular}
\end{table}

Now we test the hierarchical decomposition with $\ell$ levels. The benefit of hierarchical coarse
correction consists in the hierarchical structure and the parallelism between the coarse problems at
the same level. We fix $m_x=m_y=2$. The results are shown in Tab.~\ref{tab2}. For $\ell=1$, it is
the same as $m=2$ in Tab.~\ref{tab1}. For $\ell=2,3,4$, the subdomain number at the finest level is
the same as $m=4,8,16$ in Tab.~\ref{tab1}. Compared to Tab.~\ref{tab1}, here the coarse problems are
hierarchical, and the required $n_{c\ell}$ per subdomain at the finest level follows similar
scalings but is generally smaller. The dimension of the coarse problem for solving a subdomain
problem at level $j-1$ is half of that at level $j$. \emph{Without coarse correction and with one DD
  iteration as preconditioner}, \#iter for $\ell$ here is the same as for $m=2^\ell$ in
Tab.~\ref{tab1} because one hierarchical DD iteration is exactly one flat DD iteration in this case!

\begin{table}[t]
  \caption{GMRES \#iter \& $n_{c\ell}=2^{j-\ell}n_{cj}$ coarse modes per subdomain, $n_{i}=1$
    iteration for the coarse basis \& preconditioner, $\ell$ levels of $2\times 2$ subdomains,
    $n^24^{\ell}{\color{red}=k^2}$ elements, $4h$ overlap, {\bf free space} problem.}
  \label{tab2}
  \begin{tabular}{p{0.4cm}p{2.8cm}p{2.8cm}p{3cm}p{2.3cm}}
    \hline\noalign{\smallskip}
    $n$  & $\ell=1$ & $\ell=2$ & $\ell=3$ & $\ell=4$ \\    
    \noalign{\smallskip}\svhline\noalign{\smallskip}
    \phantom{1}4 & \phantom{1}6 \quad \phantom{1}5 \quad \phantom{1}4 \quad \phantom{1}3 & 15 \quad \phantom{1}7 \quad \phantom{1}4 \quad \phantom{1}3 & 32 \quad \phantom{1}7 \quad \phantom{1}4 \quad \phantom{1}3 & 61 \quad \phantom{1}9 \quad \phantom{1}4 \quad \phantom{1}3\\
         & \phantom{1}0 \quad \phantom{1}2 \quad \phantom{1}3 \quad \phantom{1}4 & \phantom{1}0 \quad \phantom{1}3 \quad \phantom{1}4 \quad \phantom{1}5 & \phantom{1}0 \quad \phantom{1}4 \quad \phantom{1}5 \quad \phantom{1}6 & \phantom{1}0 \quad \phantom{1}4 \quad \phantom{1}5 \quad \phantom{1}6\\
    \phantom{1}8 & \phantom{1}7 \quad \phantom{1}9 \quad \phantom{1}5 \quad \phantom{1}3 & 17 \quad \phantom{1}8 \quad \phantom{1}5 \quad \phantom{1}3 & 34 \quad \phantom{1}7 \quad \phantom{1}4 \quad \phantom{1}3 & 63 \quad 14 \quad \phantom{1}6 \quad \phantom{1}3\\
         & \phantom{1}0 \quad \phantom{1}1 \quad\phantom{1}5 \quad \phantom{1}7 & \phantom{1}0 \quad \phantom{1}6 \quad \phantom{1}7 \quad \phantom{1}8 & \phantom{1}0 \quad \phantom{1}7 \quad \phantom{1}8 \quad \phantom{1}9 & \phantom{1}0 \quad \phantom{1}7 \quad \phantom{1}8 \quad \phantom{1}9\\
    16  & \phantom{1}9 \quad \phantom{1}7 \quad \phantom{1}6 \quad \phantom{1}3 &  20 \quad \phantom{1}7 \quad \phantom{1}5 \quad \phantom{1}3 & 39 \quad \phantom{1}5 \quad \phantom{1}4 \quad \phantom{1}3 & 78 \quad \phantom{1}7 \quad \phantom{1}5 \quad \phantom{1}3  \\
         & \phantom{1}0 \quad 11 \quad 12 \quad 13 & \phantom{1}0 \quad 13 \quad 14 \quad 15 & \phantom{1}0 \quad 15 \quad 16 \quad 17 & \phantom{1}0 \quad 15 \quad 16 \quad 17\\
    32  &  10 \quad \phantom{1}5 \quad \phantom{1}4 \quad \phantom{1}3 &  22 \quad \phantom{1}6 \quad \phantom{1}4 \quad \phantom{1}3 & 48 \quad \phantom{1}5 \quad \phantom{1}4 \quad \phantom{1}3 & 92 \quad \phantom{1}5 \quad \phantom{1}4 \quad \phantom{1}3 \\
         & \phantom{1}0 \quad 27 \quad 28 \quad 29 & \phantom{1}0 \quad 29 \quad 30 \quad 31 & \phantom{1}0 \quad 31 \quad 32 \quad 33 & \phantom{1}0 \quad 32 \quad 33 \quad 34\\
    \noalign{\smallskip}\hline\noalign{\smallskip}
  \end{tabular}
\end{table}

We could fine tune $n_{cj}$ at each level and on each subdomain. For example, subdomains on the
boundary of the original domain have less interface dofs and may need less coarse modes. But we have
not done so in this paper.

\subsection{Layered Media Problem}
\label{sec:3.3}
The wavenumber in this subsection is $k(x,y)=\omega/c(y)$ piecewise constants in $y\in(0,1)$. We
assume $c(y)$ takes the values $c_1=1$ and $c_0>1$ alternatively in equally spaced layers. We keep
$\omega h=1$. The results are given in Tab.~\ref{tab3} and Tab.~\ref{tab4}.

\begin{table}
  \caption{GMRES \#iter with $n_{c}$ (below \#iter) coarse modes per subdomain, $n_{i}\!=\!1$
    iteration for the preconditioner, $m^2$ subdomains, $(mn)^2{\color{red}=\omega^2}$ elements and
    $4h$ overlap, {\bf layered media} problem.}
  \label{tab3}
  \centering $c_0=5$, 8 layers:
    \begin{tabular}{p{0.4cm}p{2.8cm}p{2.8cm}p{3cm}p{2.4cm}}
    \hline\noalign{\smallskip}
    $n$  & $m=2$ & $m=4$ & $m=8$ & $m=16$ \\    
    \noalign{\smallskip}\svhline\noalign{\smallskip}
    \phantom{1}4 & \phantom{1}7 \quad \phantom{1}5 \quad \phantom{1}4 \quad \phantom{1}3 & 19 \quad \phantom{1}6 \quad \phantom{1}4 \quad \phantom{1}3 & \phantom{1}60 \quad 10 \quad \phantom{1}5 \quad \phantom{1}3 & 103 \quad \phantom{1}6 \quad \phantom{1}4 \quad \phantom{1}3\\
         & \phantom{1}0 \quad \phantom{1}2 \quad \phantom{1}3 \quad \phantom{1}4 & \phantom{1}0 \quad \phantom{1}4 \quad \phantom{1}5 \quad \phantom{1}6 & \phantom{11}0 \quad \phantom{1}4 \quad \phantom{1}5 \quad \phantom{1}7 & \phantom{11}0 \quad \phantom{1}6 \quad \phantom{1}7 \quad \phantom{1}8\\
    \phantom{1}8 & 12 \quad \phantom{1}5 \quad \phantom{1}4 \quad \phantom{1}3 & 34 \quad \phantom{1}5 \quad \phantom{1}4 \quad \phantom{1}3 & \phantom{1}47 \quad \phantom{1}8 \quad \phantom{1}4 \quad \phantom{1}3 & 154 \quad \phantom{1}6 \quad \phantom{1}5 \quad \phantom{1}3\\
         & \phantom{1}0 \quad \phantom{1}4 \quad\phantom{1}5 \quad \phantom{1}6 & \phantom{1}0 \quad \phantom{1}9 \quad 10 \quad 11 & \phantom{11}0 \quad 10 \quad 11 \quad 13 & \phantom{11}0 \quad 11 \quad 12 \quad 13\\
    16  & 16 \quad \phantom{1}5 \quad \phantom{1}4 \quad \phantom{1}3 &  27 \quad \phantom{1}5 \quad \phantom{1}4 \quad \phantom{1}3 & \phantom{1}81 \quad \phantom{1}5 \quad \phantom{1}4 \quad \phantom{1}3 & 216 \quad \phantom{1}5 \quad \phantom{1}4 \quad \phantom{1}3  \\
         & \phantom{1}0 \quad \phantom{1}9 \quad 10 \quad 11 & \phantom{1}0 \quad 17 \quad 19 \quad 20 & \phantom{11}0 \quad 19 \quad 21 \quad 22 & \phantom{11}0 \quad 22 \quad 23 \quad 25\\
    32  &  17 \quad \phantom{1}5 \quad \phantom{1}4 \quad \phantom{1}3 &  45 \quad \phantom{1}5 \quad \phantom{1}4 \quad \phantom{1}3 & 109 \quad \phantom{1}5 \quad \phantom{1}4 \quad \phantom{1}3 & 321 \quad \phantom{1}5 \quad \phantom{1}4 \quad \phantom{1}3 \\
         & \phantom{1}0 \quad 17 \quad 20 \quad 21 & \phantom{1}0 \quad 35 \quad 37 \quad 40 & \phantom{11}0 \quad 39 \quad 40 \quad 44 & \phantom{11}0 \quad 43 \quad 44 \quad 47\\
    \noalign{\smallskip}\hline\noalign{\smallskip}
  \end{tabular}
  \smallskip
  \centering
  $c_0=10$, 8 layers:

  \begin{tabular}{p{0.4cm}p{2.8cm}p{2.8cm}p{3cm}p{2.4cm}}
    \hline\noalign{\smallskip}
    16  & 15 \quad \phantom{1}5 \quad \phantom{1}4 \quad \phantom{1}3 &  31 \quad \phantom{1}5 \quad \phantom{1}4 \quad \phantom{1}3 & \phantom{1}44 \quad \phantom{1}5 \quad \phantom{1}4 \quad \phantom{1}3 & 245 \quad \phantom{1}5 \quad \phantom{1}4 \quad \phantom{1}3  \\
         & \phantom{1}0 \quad \phantom{1}9 \quad 10 \quad 11 & \phantom{1}0 \quad 17 \quad 19 \quad 20 & \phantom{11}0 \quad 19 \quad 21 \quad 22 & \phantom{11}0 \quad 22 \quad 23 \quad 25\\
    32  &  16 \quad \phantom{1}5 \quad \phantom{1}4 \quad \phantom{1}3 &  27 \quad \phantom{1}5 \quad \phantom{1}4 \quad \phantom{1}3 & 118 \quad \phantom{1}5 \quad \phantom{1}4 \quad \phantom{1}3 & 358 \quad \phantom{1}5 \quad \phantom{1}4 \quad \phantom{1}3 \\
         & \phantom{1}0 \quad 17 \quad 20 \quad 21 & \phantom{1}0 \quad 36 \quad 37 \quad 40 & \phantom{11}0 \quad 39 \quad 40 \quad 44 & \phantom{11}0 \quad 43 \quad 44 \quad 48\\
    \noalign{\smallskip}\hline\noalign{\smallskip}
  \end{tabular}
  \smallskip
  \centering
  $c_0=10$, $64$ layers:
  
  \begin{tabular}{p{0.4cm}p{2.8cm}p{2.8cm}p{3cm}p{2.4cm}}
    \hline\noalign{\smallskip}
    32  &  14 \quad \phantom{1}7 \quad \phantom{1}5 \quad \phantom{1}3 &  58 \quad \phantom{1}5 \quad \phantom{1}4 \quad \phantom{1}3 & 235 \quad \phantom{1}5 \quad \phantom{1}4 \quad \phantom{1}3 & 185 \quad \phantom{1}5 \quad \phantom{1}4 \quad \phantom{1}3 \\
        & \phantom{1}0 \quad 21 \quad 22 \quad 23 & \phantom{1}0 \quad 33 \quad 35 \quad 38 & \phantom{11}0 \quad 39 \quad 41 \quad 43 & \phantom{11}0 \quad 40 \quad 42 \quad 44\\
    \noalign{\smallskip}\hline\noalign{\smallskip}
  \end{tabular}    
\end{table}

\begin{table}
  \caption{GMRES \#iter with $n_{c}$ coarse modes per subdomain, $n_{i}\!\!=\!\!1$ iteration for the
    coarse basis \& preconditioner, $\ell$ levels of $2\times 2$ subdomains,
    $n^24^{\ell}{\color{red}=\omega^2}$ elements, $4h$ overlap, {\bf layered media} problem.}
  \label{tab4}
  \centering
  $c_0=5$, 8 layers: 

    \begin{tabular}{p{0.4cm}p{2.8cm}p{2.8cm}p{3cm}p{2.4cm}}
    \hline\noalign{\smallskip}
    $n$  & $\ell=1$ & $\ell=2$ & $\ell=3$ & $\ell=4$ \\    
    \noalign{\smallskip}\svhline\noalign{\smallskip}
    \phantom{1}4 & \phantom{1}7 \quad \phantom{1}5 \quad \phantom{1}4 \quad \phantom{1}3 & 19 \quad \phantom{1}6 \quad \phantom{1}4 \quad \phantom{1}3 & \phantom{1}60 \quad \phantom{1}8 \quad \phantom{1}5 \quad \phantom{1}3 & 103 \quad \phantom{1}8 \quad \phantom{1}4 \quad \phantom{1}3\\
         & \phantom{1}0 \quad \phantom{1}2 \quad \phantom{1}3 \quad \phantom{1}4 & \phantom{1}0 \quad \phantom{1}3 \quad \phantom{1}4 \quad \phantom{1}5 & \phantom{11}0 \quad \phantom{1}3 \quad \phantom{1}4 \quad \phantom{1}5 & \phantom{11}0 \quad \phantom{1}4 \quad \phantom{1}5 \quad \phantom{1}6\\
    \phantom{1}8 & 12 \quad \phantom{1}5 \quad \phantom{1}4 \quad \phantom{1}3 & 34 \quad \phantom{1}6 \quad \phantom{1}4 \quad \phantom{1}3 & \phantom{1}47 \quad \phantom{1}6 \quad \phantom{1}4 \quad \phantom{1}3 & 154 \quad \phantom{1}7 \quad \phantom{1}5 \quad \phantom{1}3\\
         & \phantom{1}0 \quad \phantom{1}4 \quad\phantom{1}5 \quad \phantom{1}6 & \phantom{1}0 \quad \phantom{1}5 \quad \phantom{1}6 \quad \phantom{1}7 & \phantom{11}0 \quad \phantom{1}6 \quad \phantom{1}7 \quad \phantom{1}8 & \phantom{11}0 \quad \phantom{1}7 \quad \phantom{1}8 \quad \phantom{1}9\\
    16  & 16 \quad \phantom{1}5 \quad \phantom{1}4 \quad \phantom{1}3 &  27 \quad \phantom{1}7 \quad \phantom{1}5 \quad \phantom{1}3 & \phantom{1}81 \quad \phantom{1}7 \quad \phantom{1}4 \quad \phantom{1}3 & 216 \quad 10 \quad \phantom{1}5 \quad \phantom{1}3  \\
         & \phantom{1}0 \quad \phantom{1}9 \quad 10 \quad 11 & \phantom{1}0 \quad 10 \quad 11 \quad 12 & \phantom{11}0 \quad 12 \quad 13 \quad 14 & \phantom{11}0 \quad 14 \quad 15 \quad 16\\
    32  &  17 \quad \phantom{1}5 \quad \phantom{1}4 \quad \phantom{1}3 &  45 \quad \phantom{1}5 \quad \phantom{1}4 \quad \phantom{1}3 & 109 \quad \phantom{1}5 \quad \phantom{1}4 \quad \phantom{1}3 & 321 \quad \phantom{1}6 \quad \phantom{1}5 \quad \phantom{1}3 \\
         & \phantom{1}0 \quad 17 \quad 20 \quad 21 & \phantom{1}0 \quad 23 \quad 24 \quad 25 & \phantom{11}0 \quad 26 \quad 27 \quad 28 & \phantom{11}0 \quad 28 \quad 29 \quad 30\\
    \noalign{\smallskip}\hline\noalign{\smallskip}
  \end{tabular}
\end{table}


\bibliographystyle{spmpsci}
\bibliography{ms01_zhang.bib}

\end{document}